\documentclass{amsart}
\usepackage{eurosym}
\usepackage{amsfonts}
\usepackage{amsmath}
\usepackage{amssymb}
\usepackage{graphicx}
\usepackage{amscd}

\newtheorem{theorem}{Theorem}
\theoremstyle{plain}

\newtheorem{definition}{Definition}

\newtheorem{lemma}{Lemma}

\newtheorem{remark}{Remark}

\numberwithin{equation}{section}

\begin{document}
\title[General Central Limit Theorems]{General Central Limit Theorems for
Associated Sequences}
\author{$^{(1,2)}$ Harouna SANGARE}
\author{$^{(1,3,4)}$ Gane Samb LO}
\email[Harouna SANGARE]{harounasangareusttb@gmail.com}
\email[Gane Samb LO]{ganesamblo@ganesamblo.net}

\begin{abstract}
In this paper, we provide general central limit theorems (\textit{CLT}'s) for associated random variables (\textit{rv}'s) following the approaches used by Newman (1980) and Olivera et al.(2012). Given some assumptions, a Lyapounov-Feller-Levy type theorem is stated. We next specify different particular \textit{CLT} versions of associated sequences based on moment conditions. A comparison study with available \textit{CTL}'s is performed. As a by-product, we complete an important available theorem where an assumption was missing.\\

\bigskip \noindent $^{(1)}$ LERSTAD, Gaston Berger University, Saint-Louis,
Senegal.\newline
\noindent $^{(2)}$ DER de Math\'ematiques et d'Informatique, USTTB, Bamako,
Mali.\newline
\noindent $^{(3)}$ AUST - African University of Sciences and Technology,
Abuja, Nigeria\newline
$^{(4)}$ LSTA, Pierre and Marie Curie University, Paris VI, France.\newline

\noindent \textit{Corresponding author}. Gane Samb Lo. Email :
gane-samb.lo@edu.ugb.sn. Postal address : 1178 Evanston Dr NW T3P 0J9,
Calgary, Alberta, Canada.
\end{abstract}

\maketitle

\Large

\section{Introduction} \label{sec1}

\noindent We consider the problem of the central limit theorem for
associated sequences. This problem goes back to Newman \cite{newman80}.
Since then, a number of $CLT$'s, and strong laws of large numbers (\textit{SLLN}%
's) or weak laws of large numbers (\textit{WLLN}'s), and invariance
principles and laws of the iterated logarithm (\textit{LIL}'s) have been
provided in the recent literature by different authors. Dabrowski and
co-authors (see \cite{burton} and \cite{dabro}) considered weakly associated
random variables to establish invariance principles in the lines of Newman
and Wright \cite{newmanwright}, as well as Berry-Essen-type results and
functional \textit{LIL}'s. The weak convergence for empirical processes
using associated sequences has been discussed by Louhichi \cite{louhichi}
and Yu \cite{yu93}.\newline

\noindent The most general \textit{CLT}'s seem to be the one provided by Cox
and Grimmet \cite{cox} for arbitrary associated \textit{rv}'s fulfilling a number of
moment conditions and those given by Oliveira \cite{paulo}.\newline

\noindent This question arises in the active research field on the concept
of association and its application in many sciences, especially in
percolation theory in Physics and in Reliability. The books by Rao \cite%
{rao} and the monograph by Oliveira \cite{paulo} present a review of these
researches. The book of Bulinski and Shashkin \cite{bulinski2007} treats
random associated sequences and intensively uses properties of regularly
varying functions and provides \textit{CTL}'s, \textit{LLN}'s, \textit{LIL}%
's and invariance principles.\newline

\noindent Although many results concerning the \textit{CLT} problem are
available for such sequences, there are still a number of open problems,
especially regarding nonstationary sequences. \bigskip

\noindent Here, we intend to provide more general \textit{CLT}'s for arbitrary associated
sequences. Precisely, here, we want to use all the power of the Newman
method and express the conditions in the most general frame based on
moment conditions so that any other result might be derived from them. In such a
way, a Lyapounov-Feller-Levy type Theorem will be possible to be stated given
some general assumptions. From this approach, more general \textit{CLT}'s
may be conceived only by turning back the Newman's method.\\

\noindent The paper is organized as follows. Since association is the
central notion used here, we first make a quick reminder of it in Section %
\ref{sec2}. In Section \ref{sec3}, we make a round up of CLT's available in
the literature with the aim of comparing them to our findings. In Section
\ref{sec4}, we state our general \textit{CLT} versions for arbitrary associated \textit{rv}'s and give 
comparison study and concluding remarks.

\bigskip

\section{A brief Reminder of Association} \label{sec2}

\label{sec2}\noindent The concept of association has been introduced by
Lehmann (1966) \cite{lehmann} in the bivariate case, and extended later in
the multivariate case by Esary, Proschan and Walkup (1967) \cite{esary}.
\noindent The concept of association for random variables generalizes that
of independence and seems to model a great variety of stochastic models.%
\newline

\noindent This property also arises in Physics, and is quoted under the name
of FKG property (Fortuin, Kastelyn and Ginibre (1971) \cite{fortuin}), in
percolation theory and even in Finance (see \cite{jiazhu}).\newline

\noindent The definite definition is given by Esary, Proschan and Walkup
(1967) \cite{esary} as follows.

\begin{definition}
A finite sequence of rv's $(X_{1},...,X_{n})$ is associated if, for any
couple of real and coordinate-wise non-decreasing functions $h$ and $g$
defined on $\mathbb{R}^{n}$, we have 
\begin{equation}
Cov(h(X_{1},...,X_{n}),\ \ g(X_{1},...,X_{n}))\geq 0,  \label{asso}
\end{equation}

\noindent whenever the covariance exists. An infinite sequence of rv's are
associated whenever all its finite subsequences are associated.
\end{definition}

\noindent We have a few number of interesting properties to be found in (\cite{rao}) :\\

\noindent \textbf{(P1)} A sequence of independent rv's is associated.\newline

\noindent \textbf{(P2)} Partial sums of associated rv's are associated.%
\newline

\noindent \textbf{(P3)} Order statistics of independent rv's are associated.%
\newline

\noindent \textbf{(P4)} Non-decreasing functions and non-increasing
functions of associated variables are associated.\newline

\noindent \textbf{(P5)} Let the sequence $Z_{1},Z_{2},...,Z_{n}$ be
associated and let $(a_{i})_{1\leq i\leq n}$ be positive numbers and $%
(b_{i})_{1\leq i\leq n}$ real numbers. Then the \textit{rv}'s $%
a_{i}(Z_{i}-b_{i})$ are associated.\newline

\noindent As immediate other examples of associated sequences, we may cite
Gaussian random vectors with nonnegatively correlated components (see \cite%
{pitt}) and a homogeneous Markov chain (see \cite{daley}).\newline

\noindent Demimartingales are set from associated centered variables exactly
as martingales are derived from partial sums of centered independent random
variables. We have

\begin{definition}
A sequence of rv's $\{S_{n},n\geq 1\}$ \ in $L^{1}(\Omega ,\mathcal{A},%
\mathbb{P})$ \ is a demimartingale when for any $j\geq 1$, for any
coordinatewise nondecreasing function $g$ \ defined on $\mathbb{R}^{j}$, we
have 
\begin{equation}
\mathbb{E}\left( (S_{j+1}-S_{j})\ g(S_{1},...,S_{j})\right) \geq 0,\ \ j\geq
1.  \label{defmarting}
\end{equation}
\end{definition}

\bigskip \noindent Two particular cases should be highlighted. First any
martingale is a demimartingale. Secondly, partial sums $S_{0}=0$, $%
S_{n}=X_{1}+...+X_{n}$, $n\geq 1$, of associated and centered random
variables $X_{1},X_{2},...$ are demimartingales. In this case, (\ref%
{defmarting}) becomes :

\begin{equation*}
\mathbb{E}\left( (S_{j+1}-S_{j})\ g(S_{1},...,S_{j})\right) =\mathbb{E}%
\left( X_{j+1}\ g(S_{1},...,S_{j})\right) =Cov\left(
X_{j+1},g(S_{1},...,S_{j})\right) ,
\end{equation*}

\noindent since $\mathbb{E}X_{j+1}=0$. Since $(x_{1},...,x_{j+1})\longmapsto
x_{j+1}$ and $(x_{1},...,x_{j+1})\longmapsto g(x_{1},...,x_{j})$ are
coordinate-wise nondecreasing functions and since the $X_{1},X_{2},..$ are
associated, we get

\begin{equation*}
\mathbb{E}\left( (S_{j+1}-S_{j})\ g(S_{1},...,S_{j})\right) =Cov\left(
X_{j+1}\ g(S_{1},...,S_{j})\right) \geq 0.
\end{equation*}

\bigskip \noindent Finally, we present the following key results for
associated sequences that one can find in almost any paper on that topic and
that we need for our proofs. A detailed\ review on these results is given in 
\cite{gslo}.

\begin{lemma}[Hoeffding (1940) see \cite{rao}]
\label{lemg1} Let $(X,Y)$ be a bivariate random vector such that $\mathbb{E}%
(X^{2})<+\infty $ and $\mathbb{E}(Y^{2})<+\infty .$ If $\left(
X_{1},Y_{1}\right) $ and $\left( X_{2},Y_{2}\right) $ are two independent
copies of $(X,Y),$ then we have

\begin{equation*}
2Cov(X,Y)=\mathbb{E}(X_{1}-X_{2})(Y_{1}-Y_{2}).
\end{equation*}

\noindent We also have 
\begin{equation*}
Cov(X,Y)=\int_{-\infty }^{+\infty }\int_{-\infty }^{+\infty }H(x,y)dxdy,
\end{equation*}
\end{lemma}

\noindent where 
\begin{equation*}
H(x,y)=\mathbb{P}(X>x,Y>y)-\mathbb{P}(X>x)\mathbb{P}(Y>y).
\end{equation*}

\bigskip

\begin{lemma}[Newman (1980), see \cite{newman80}]
\label{lemg2} Suppose that $X$, $Y$ are two random variables with finite
variance and, $f$ and $g$ are $\mathbb{C} ^{1}$ complex valued functions on $%
\mathbb{R}^{1}$ with bounded derivatives $f^{\prime }$ and $g^{\prime }.$
Then

\begin{equation*}
|Cov(f(X),g(Y))|\leq ||f^{\prime }||_{\infty }||g^{\prime }||_{\infty
}Cov(X,Y).
\end{equation*}
\end{lemma}

\noindent The following lemma is the most used tool in this field.

\begin{lemma}[Newman and Wright (1981) Theorem, see \protect\cite%
{newmanwright}]
\label{lemg3} Let $X_{1},X_{2},...,X_{n}$ be associated, then we have for
all $t=(t_{1},...,t_{n})\in \mathbb{R}^{k}$,

\begin{equation}
\left\vert \psi _{_{(X_{1},X_{2},...,X_{n})}}(t)-\prod\limits_{i=1}^{n}\psi
_{_{X_{i}}}(t_{i})\right\vert \leq \frac{1}{2}\sum_{1\leq i\neq j\leq
n}\left\vert t_{i}t_{j}\right\vert Cov(X_{i},X_{j}).  \label{decomp}
\end{equation}
\end{lemma}

\bigskip \noindent Before we proceed any further, let us  make a round up of \textit{%
CLT}'s for associated sequences in stationary and non-stationary cases.

\section{Central limit theorem for associated sequences} \label{sec3} \noindent Let $X_{1},X_{2},\cdots ,X_{n}$ be an associated
sequence of mean-zero random variables defined on the same probability space
($\Omega ,\mathcal{A},\mathbb{P}$). Define for each $n\geq 1,$ $%
S_{n}=X_{1}+...+X_{n}.$ The CLT question for stationary associated sequence
turns around Newman (see \cite{newman80}) results in which is proved that $%
S_{n}/\sqrt{n}$ converges to a normal random variable $\mathcal{N}(0,\sigma
^{2})$ when

\begin{equation*}
\sigma ^{2}=\mathbb{V}ar(X_{1})+2\sum_{j=2}^{\infty
}Cov(X_{1},X_{j})<+\infty .
\end{equation*}

\noindent And in such a situation, 
\begin{equation*}
s_{n}^{2}=\mathbb{V}ar(S_{n}/\sqrt{n})\rightarrow \sigma ^{2}\text{ as }%
n\rightarrow +\infty .
\end{equation*}

\noindent A number of invariance principles and other CLT's are available
but they are generally adaptations of this Newman result. As to the general
case, Cox and Grimmet (see \cite{cox}), did not consider stationarity in
their results which used triangular sequences. Formulated for simple
sequences, their result is that $S_{n}/(s_{n}\sqrt{n})$ weakly converges to a normal
random variable $\mathcal{N}(0,1)$ if $\ \mathbb{V}ar(X_{n})$ is
asymptotically bounded below zero, and the sequence of the third moments $%
\mathbb{E}\left\vert X_{n}\right\vert ^{3}$ is bounded and there exists a
function $u(r)$, $r\in \{0,1,...\}$ such that \ $u(r)\rightarrow 0$ as $%
r\rightarrow +\infty $ and such for all $k\geq 1$, and all $r\geq 0$

\begin{equation*}
\sum_{j:|k-j|\geq r}Cov(X_{j},X_{k})\leq u(r).
\end{equation*}

\noindent Let us recall his $CLT$ as follows

\bigskip

\begin{theorem}
\label{cox} Let $X_{1},X_{2},\cdots ,X_{n}$ be an associated sequence of
mean-zero random variables defined on the same probability space ($\Omega ,%
\mathcal{A},\mathbb{P}$). Suppose there exist positive and finite constants $%
c_{1}$ and $c_{2}$ such that

\begin{equation}
Var(X_{j})\geq c_{1}\text{ and }\mathbb{E}\left\vert X_{j}\right\vert
^{3}\leq c_{2}\text{ for all }j\geq 1  \label{coxA1},
\end{equation}

\noindent and there is a function $u(r)$ of $r\in \mathbb{N}$ such that $%
u(r)\rightarrow 0$ as r$\rightarrow +\infty $ and for any $r\geq 1,$%
\begin{equation}
\sup_{j\geq 1}\sum_{i:\left\vert j-i\right\vert \geq r}cov(X_{i},X_{j})\leq
u(r).  \label{coxA2}
\end{equation}

\noindent Then 
\begin{equation*}
S_{n}/s_{n}\rightsquigarrow N(0,1)\text{ as }n\rightarrow +\infty .
\end{equation*}

\noindent where throughout the text, the symbol $\rightsquigarrow $ stands
for the weak convergence.
\end{theorem}

\bigskip\ 

\noindent Oliveira \textsl{et al.} \cite{paulo} have proved general \textit{CLT}'s, still using
the Newman approach.\\

\noindent First, they obtained 

\begin{theorem}[see \protect\cite{paulo}, page 105, Theorem 4.4]
\label{theooliv1} Let $X_{n}$, $n\in \mathbb{N}$, be centered,
square-integrable and associated random variables. For each $n\in \mathbb{N}$%
, let $\ell _{n}\in \mathbb{N}$ and $m_{n}=$ $\left[ \frac{n}{\ell _{n}}%
\right] $. Define, for $j=1,...,m_{n},$ $Y_{j,\ell _{n}}=\sum_{i=\left(
j-1\right) \ell _{n}+1}^{j\ell _{n}}X_{i}$ and $Y_{m_{n}+1,\ell
_{n}}=\sum_{i=m_{n}\ell _{n}+1}^{n}X_{i}$. Assume that $m_{n}\rightarrow
+\infty$, and that

\begin{equation}
\frac{1}{s_{n}^{2}}\sum_{j=1}^{m_{n}}\mathbb{V}ar\left( Y_{j,\ell
_{n}}\right) \rightarrow 1,  \label{paulaA1}
\end{equation}
\begin{equation}
\left\vert \mathbb{E}\exp \left( \frac{it}{s_{n}}S_{n}\right)
-\prod\limits_{j=1}^{m_{n}}\mathbb{E}\exp \left( \frac{it}{s_{n}}Y_{j,\ell
_{n}}\right) \right\vert \rightarrow 0,\text{ }t\in \mathbb{R},
\label{pauloA2}
\end{equation}

\noindent and

\begin{equation}
\forall \text{ }\varepsilon >0,\text{ }\frac{1}{s_{n}^{2}}%
\sum_{j=1}^{m_{n}}\int_{\left\{ \left\vert Y_{j,\ell _{n}}\right\vert \geq
\varepsilon s_{n}\right\} }Y_{j,\ell _{n}}d\mathbb{P}\rightarrow 0.
\label{pauloA3}
\end{equation}

\noindent Then 
\begin{equation*}
\frac{1}{s_{n}}S_{n}\rightsquigarrow \mathcal{N}(0,1),
\end{equation*}
\end{theorem}

\noindent Next, they obtained the following result using a Feller-Levy
condition.

\begin{theorem}[see \protect\cite{paulo}, page 108, Theorem 4.8]
\label{theooliv2} Let $X_{n}$, $n\in \mathbb{N}$, be centered,
square-integrable and associated random variables. Assume that

\begin{equation}
u\left( n\right) \rightarrow 0,\text{ }u\left( 1\right) <+\infty ,
\label{pauloB1}
\end{equation}%
\begin{equation}
\inf_{n\in N}\frac{1}{n}s_{n}^{2}>0,  \label{pauloB2}
\end{equation}

\begin{equation}
\forall \text{ }\varepsilon >0,\text{ }\frac{1}{s_{n}^{2}}%
\sum_{j=1}^{m_{n}}\int_{\left\{ \left\vert X_{j}\right\vert \geq \varepsilon
s_{n}\right\} }X_{j}^{2}d\mathbb{P}\rightarrow 0.  \label{pauloB3}
\end{equation}

\noindent Then 
\begin{equation*}
\frac{1}{s_{n}}S_{n}\rightsquigarrow \mathcal{N}(0,1).
\end{equation*}
\end{theorem}

\noindent \textbf{Remark on whether the assumptions of the theorem are enough to get the \textit{CLT}}. It seems to us that the conditions given by this theorem are not enough, as we tried to show it in Subsection \ref{subsec42} below. We think that the following assumption, denoted (\textit{Hab}) below,

\begin{equation*}
\frac{1}{s_{n}^{2}}\mathbb{V}ar\left( \sum_{i=m(n)\ell
(n)+1}^{r(n)}X_{j}\right) \rightarrow 0\text{ as }n\rightarrow +\infty .
\end{equation*}

\noindent should be added. This assumption is implied by this simpler one, denoted (\textit{HNab}) below : 
\begin{equation*}
\frac{1}{s_{n}^{2}}\sum_{i=t_{n}}^{u_{n}}var(X_{i})\rightarrow 0\text{ as }
n\rightarrow \infty ,
\end{equation*}

\noindent for $0\leq t_{n}\leq u_{n}\leq n,u_{n}-t_{n}\leq \ell(n),(u_{n}-v_{n})/n\rightarrow 0$ as $n\rightarrow \infty$.\\

\noindent For a stationary case, this assumption is immediate. The foundation of our remark is given in Point (1) in Subsubsection 
\ref{subsubsec423} of Subsection \ref{subsec42} in Section \ref{sec4}.\\ 

\bigskip \noindent Our objective in this paper is to express \textit{CLT's} in the most general
setting, still using the Newman approach and to derive the former results as particular cases. With respect to the former results described above, we simplify the approach and get the best we can do by formulating a Lyapounov-Feller-Levy type of \textit{CLT}. The general conditions are next expressed on moment conditions stated also in a general setting. Existing versions are all included in our statements. And we conclude that more general \textsl{CLT}'s cannot be obtained without getting out the Newman approach.\\

\noindent Let us begin by introduce the following assumptions

\bigskip

\noindent There exists a sequence $\ell (n)$ of positive integers such that $%
n=m(n)\ell (n)+r(n),$ with $0\leq r(n)<\ell (n),$ $0\leq m(n)\rightarrow
+\infty $ and 
\begin{equation}
(\ell (n)/n,r(n)/n)\rightarrow (0,0)\text{ \ }as\text{ \ }n\rightarrow
+\infty .  \tag{L}
\end{equation}

\bigskip \noindent We want to stress that the integers $m=m(n),$ $\ell =\ell (n)$ and $r=r(n)$
depend of $n$ throughout the text even though we may and do drop the label $n
$ in many situations for simplicty's sake.\\

\noindent On top of this general assumption, we may require the following
ones.

\begin{equation}
\frac{\ell (n)}{s_{n^{2}}}\rightarrow 0\text{ as }n\rightarrow +\infty .
\label{H0}
\end{equation}%
\begin{equation}
\frac{\ell (n)}{s_{n}^{2}}\sum_{j=1}^{m(n)}\mathbb{V}ar\left( \frac{S_{j\ell
(n)}-S_{(j-1)\ell (n)}}{\sqrt{\ell (n)}}\right) \rightarrow 1\text{ as }%
n\rightarrow +\infty .  \tag{Ha}
\end{equation}

\begin{equation}
\frac{1}{s_{n}^{2}}\mathbb{V}ar\left( \sum_{i=m(n)\ell
(n)+1}^{r(n)}X_{j}\right) \rightarrow 0\text{ as }n\rightarrow +\infty . 
\tag{Hab}
\end{equation}

\noindent 

\begin{equation}
\sup_{1\leq j\leq m(n)+1}\frac{\ell (n)}{s_{n}^{2}}\mathbb{V}ar\left( \frac{%
S_{j\ell (n)}-S_{(j-1)\ell (n)}}{\sqrt{\ell (n)}}\right)
=C_{1}(n)\rightarrow 0\text{ as \ }j\rightarrow +\infty .  \tag{Hb}
\end{equation}

\noindent We have for some $\delta >0,$  $\mathbb{E}\left\vert
X_{j}\right\vert ^{2+\delta }<+\infty ,j\geq 1$ and the Lyapounov Condition
holds

\begin{equation}
\frac{\ell ^{3/2}(n)}{s_{n}^{2+\delta }}\sum_{j=1}^{m}\mathbb{E}\left\vert 
\frac{S_{j\ell (n)}-S_{(j-1)\ell (n)}}{\sqrt{\ell (n)}}\right\vert
^{2+\delta }=C_{2}(n)\rightarrow 0\text{ as j}\rightarrow +\infty .  \tag{Hc}
\end{equation}

\noindent In the sequel, it may be handy to use the notation 
\begin{equation}
Y_{j,\ell }=\frac{S_{j\ell (n)}-S_{(j-1)\ell (n)}}{\sqrt{\ell (n)}},1\leq
j\leq m=m(n).  \label{THEY}
\end{equation}

\bigskip

\noindent In the next section, we will state two \textit{CLT}'s based on these hypotheses and next. A third ine will be a completition of Theorem \ref{theooliv2}. Next, the results are particularized into more specific versions.\\

\section{Results and Commentaries} \label{sec4} \noindent In this section, we present general \textit{CLT}'s for
associated \textit{rv}'s and next give different forms in specific types of independent and dependent data and
finally make a comparison with available results.

\subsection{General \textit{CLT}'s} \label{subsec41} \noindent We have following results.

\begin{theorem}
\label{theo1} Let $X_{1},X_{2},\cdots ,X_{n}$ be an associated sequence of
mean-zero random variables defined on the same probability space ($\Omega ,%
\mathcal{A},\mathbb{P}$). If the sequence is stationary, then 
\begin{equation*}
\dfrac{S_{n}}{\sqrt{n}}=\dfrac{X_{1}+X_{2}+\cdots +X_{n}}{\sqrt{n}}%
\rightsquigarrow \mathcal{N}(0,\sigma )\ as\ n\rightarrow +\infty ,
\end{equation*}

\noindent In the general setting, if \textit{(L)}, (H0), \textit{(Ha)}, 
\textit{(Hb)} and \textit{(Hc)} hold, then 
\begin{equation*}
\dfrac{S_{n}}{s_{n}}=\dfrac{X_{1}+X_{2}+\cdots +X_{n}}{s_{n}} \rightsquigarrow \mathcal{N}(0,1)\ as\ n\rightarrow +\infty .
\end{equation*}
\end{theorem}

\bigskip

\bigskip \noindent Next, we state a Lyapounov-Feller-Levy type theorem given
some assumptions.

\begin{theorem}
\label{theo2} Let $X_{1},X_{2},\cdots ,X_{n}$ be an associated sequence of
mean-zero random variables defined on the same probability space ($\Omega ,%
\mathcal{A},\mathbb{P}$). Denote for each $j\in \{1,m\},\tau
_{j}^{2}=Var\left( S_{j\ell }-S_{(j-1)\ell }\right) =\mathbb{E}\left(
S_{j\ell }-S_{(j-1)\ell }\right) ^{2}$ and%

\begin{equation*}
\nu _{m(n)}^{2}=\tau _{1}+...+\tau _{m(n)},n\geq 1\text{.}
\end{equation*}

\noindent Assume that the assumptions \textit{(L)} and \textit{(Ha) hold and }
either (Hab) or \textit{(Hb) is true. The }we have the following equivalence
result :

\begin{equation*}
\max_{1\leq k\leq m(n)}^{2}\mathbb{E}\left( S_{j\ell }-S_{(j-1)\ell }\right)
/s_{n}^{2}\rightarrow 0\text{ as n}\rightarrow +\infty
\end{equation*}

\noindent and 
\begin{equation*}
S_{n}/s_{n}\rightsquigarrow \mathcal{N}(0,1)\ as\ n\rightarrow +\infty ,
\end{equation*}

\noindent if and only if for any $\varepsilon >0$,

\begin{equation}
\frac{1}{s_{n}^{2}}\mathbb{E}\left( \left( S_{j\ell }-S_{(j-1)\ell }\right)
^{2}1_{\left( \left\vert S_{j\ell }-S_{(j-1)\ell }\right\vert \geq
\varepsilon \nu _{m(n)}\right) }\right) \rightarrow 0\text{ as }n\rightarrow
0.  \label{fellerLevy01b}
\end{equation}

\noindent Moreover the sequence $(\nu _{m(n)})_{n\geq 1}$ may be replaced by the sequence of 
$(s_{n})_{n\geq 1}$\ in Condition (\ref{fellerLevy01b}) to give

\begin{equation}
\frac{1}{s_{n}^{2}}\mathbb{E}\left( \left( S_{j\ell }-S_{(j-1)\ell }\right)
^{2}1_{\left\vert S_{j\ell }-S_{(j-1)\ell }\right\vert \geq \varepsilon
s_{n}}\right) \rightarrow 0\text{ as }n\rightarrow 0.  \label{fellerlevy01c}
\end{equation}
\end{theorem}

\begin{theorem}
\label{theo3} Let $X_{1},X_{2},\cdots ,X_{n}$ be an associated sequence of
mean-zero random variables defined on the same probability space ($\Omega ,%
\mathcal{A},\mathbb{P}$). If 
\begin{equation*}
\left\vert \Psi _{\frac{S_{m\ell }}{s_{n}}}(t)-\prod\limits_{j=1}^{m}\Psi
_{Y_{j,n}}\left( \frac{\sqrt{\ell }}{s_{n}}t\right) \right\vert \rightarrow 0%
\text{ as }n\rightarrow +\infty,
\end{equation*}

\noindent then we have the following equivalence result :

\begin{equation*}
\max_{1\leq k\leq m(n)}^{2}\mathbb{E}\left( S_{j\ell }-S_{(j-1)\ell }\right)
/s_{n}^{2}\rightarrow 0\text{ as n}\rightarrow +\infty
\end{equation*}

\noindent and 
\begin{equation*}
S_{n}/s_{n}\rightsquigarrow \mathcal{N}(0,1)\ as\ n\rightarrow +\infty ,
\end{equation*}

\noindent if and only if for any $\varepsilon >0$, \ Formula \ref%
{fellerLevy01b} holds.
\end{theorem}

\bigskip 

\begin{remark} \noindent Let us make the following remarks.\newline

\noindent \textbf{(1)} The method of proving the aboved theorem will consist
of decomposing the sums of variables into sums of blocks of variables and
treating these as if they were independent. Naturally, we will need some
control on the approximation between the sums of the dependent blocks and
their independent counterparts. This control will be achieved using
characteristic functions and is based on the inequality in Lemma \ref{lemg3} of Newman. Our approach is to go the farest possible only using the
moments conditions.\newline

\noindent \textbf{(2)} Theorem \ref{theo2} is not yet a
Lyapounov-Feller-Levy Theorem $(LFLT)$. Using Lemma \ref{lemg3}, its only
say we have a $LFLT$ provided assumptions that make the CLT problem into a
CLT one concerning indedendent variables. A full $LFLT$ can not be achieved
as long as the proofs are based on the approximation of Lemma \ref{lemg3}.
\end{remark}

\noindent Before we proceed to the proofs in Subsection \ref{subsec5}, we
are now going to derive some consequences and particular cases of the
theorem.

\subsection{Commentaries and Consequences} \label{subsec42}

\subsubsection{The most general approach leading to a Feller-Levy CLT type} \label{subsubsec421} We begin a general comments of the approach.\\

\noindent Almost all the available \textit{CLT} results used the Newman's method based on Lemma \ref{lemg3}. The approach we used is intended to get the sharpest results we can get in that frame. In ealier versions of our results, we were not aware of the
results of Oliveira \textit{et al. \cite{paulo}.} However, with the
knowledge of these results, our works still present a number of significant
avantages we want to highlight here. Actually, Oliveira et al. \cite{paulo}
attains the best we can do in the Newman approach : the best and unique way to find out assumptions under which
\begin{equation*}
S_{n}/s_{n}\rightarrow N(0,1),
\end{equation*}

\noindent holds, is to reduce it to
\begin{equation}
\prod\limits_{j=1}^{m}\Psi _{Y_{j,n}}\left( \frac{\sqrt{\ell }}{s_{n}}%
t\right) \rightarrow \exp (-t^{2}/2)\text{ }as\text{ }n\rightarrow +\infty .
\label{C1}
\end{equation}

\noindent after are defined the random variables $Y_{j,\ell }$ of characteristic functions $\Psi _{Y_{j,n}}$ given in Formula (\ref{THEY}) based on the decomposition \textit{(L)} (We recall that our notation of $Y_{j,\ell }$ are not the sams as that of  \cite{paulo}).\\ 

\noindent This is the justification of Assumption (\ref{pauloA2}) above, which
corresponds to the equivalent one we used, which is 
\begin{equation}
\left\vert \Psi _{\frac{S_{m\ell }}{s_{n}}}(t)-\prod\limits_{j=1}^{m}\Psi
_{Y_{j,n}}\left( \frac{\sqrt{\ell }}{s_{n}}t\right) \right\vert \rightarrow 0%
\text{ as }n\rightarrow +\infty .  \label{C2}
\end{equation}

\noindent From there, the authors did not, as far as we know, capitalize this fact in order to have a Feller-Levy final \textit{CLT} version, as we did in Theorem \ref{theo3}. In our view, this version is the starting point for new \textit{CLT}'s out of the Newman
approach.

\subsubsection{General condition} \label{subsubsec422}

Based on the best the Newman approach can give, it remains to have the most
general conditions to ensure (\ref{C1}) and (\ref{C2}). If we wish to
directly express (\ref{C2}) into the  $X_{i}$'s, Oliveira \textit{et al.} 
\cite{paulo} proved in page 109 that their assumption (\ref{pauloB3})
implies Formula (\ref{fellerlevy01c}) in Theorem \ref{theo2} above. In
general, authors usually provide \textit{CTL}'s based on conditions ensuring 
(\ref{C1}) and (\ref{C2}).\\ 

\noindent In that specific case, we proceeded into two directions :\\

\noindent \textbf{(1)} Expressing general conditions based on moments. We will see in the next
subsection, how the available \textit{CLT}'s may be derived from Theorem \ref{theo2}.\\

\noindent \textit{(2)} Keeping the notation of decomposition $(L)$ in the assumptions. This
will allow, in particular cases, to base methods on specific values of $m(n)$ and $\ell (n)$ 

\subsubsection{Comparisons} \label{subsubsec423} Let us highlight some comparison results.\\

\noindent \textbf{(1)} \textbf{With Theorem \ref{theooliv2}. A possible gap in Theorems \ref{theooliv2} of \cite{paulo}}.\\

\noindent By setting

\begin{equation*}
Y_{j,\ell (n)}^{\ast }=Y_{j,\ell (n)}\text{ for }j=1,...,m(n)\text{ and }%
Y_{m+1,\ell }^{\ast }=\sum_{i=m(n)\ell (n)+1}^{r(n)}X_{i}/\sqrt{\ell (n)},
\end{equation*} 

\noindent we have  
\begin{equation*}
S_{n}=\ell (n)\sum_{j=1}^{m(n)+1}Y_{j,\ell (n)}^{\ast }.
\end{equation*}

\noindent It comes that

\begin{eqnarray*}
s_{n}^{2}=Var(S_{n})&=&\ell (n)\sum_{j=1}^{m(n)}var(Y_{j,\ell (n)})+\ell
^{2}Var(Y_{m(n)+1,\ell (n)}^{\ast })\\
&+&2\ell (n)\sum_{1\leq h<k\leq m(n)+1}cov(Y_{h,\ell (n)}^{\ast },Y_{k,\ell (n)}^{\ast }).
\end{eqnarray*}

\noindent By definition of the Cox coefficient (as denamed by Bulinski \textit{et al.}
\ \cite{bulinski2007}\ ), and since the indices of the $X_{i}$'s in $%
cov(Y_{h,\ell }^{\ast },Y_{k,\ell }^{\ast })$ are distanced by $1,...,\ell $
points in absolute values, and the $X_{i}$'s are therein normalized by $%
\sqrt{\ell }$, we have
\begin{eqnarray*}
\sum_{1\leq h<k\leq m(n)+1}cov(Y_{h,\ell (n)}^{\ast },Y_{k,\ell (n)}^{\ast
}) &\leq &\frac{1}{\ell }\sum_{h=1}^{\ell (n)}\sup_{i\geq
1}\sum_{k:\left\vert k-i\right\vert \geq h}cov(X_{i},X_{k}) \\
&\leq &\frac{1}{\ell} \sum_{i=1}^{\ell(n)}u(i),
\end{eqnarray*}

\noindent and then

\begin{equation*}
s_{n}^{2}=Var(S_{n})\leq \ell \sum_{j=1}^{m(n)}var(Y_{j,\ell (n)})+\ell
(n)^{2}Var(Y_{m(n)+1,\ell (n)}^{\ast })+2\sum_{i=1}^{\ell (n)}u(i),
\end{equation*}

\noindent which gives
\begin{eqnarray*}
\left\vert 1-\frac{\ell (n)}{s_{n}^{2}}\sum_{j=1}^{m}var(Y_{j,\ell (n)}-%
\frac{\ell (n)^{2}}{s_{n}^{2}}Var(Y_{m(n)+1,\ell (n)}^{\ast })\right\vert 
&\leq &\frac{2\ell (n)}{s_{n^{2}}}\left\{ \frac{1}{\ell (n)}\sum_{i=1}^{\ell(n)}u(i)\right\}  \\
&\rightarrow &0\text{ as }n\rightarrow \infty ,
\end{eqnarray*}

\noindent by C\'{e}saro's Lemma if 
\begin{equation}
\lim \sup_{n\rightarrow +\infty }\frac{\ell (n)}{s_{n^{2}}}<+\infty .
\label{pauloB2S}
\end{equation}

\noindent Here, it seems to us that the authors of \cite{paulo} might have not taken
into a account the term $\ell Y_{m+1,\ell }^{\ast }$ in the line --7 of
their page 108. At line -6 of the same page, their  formula $S_{n}=Y_{1,\ell
_{n}}+...+Y_{m_{n},\ell _{n}}$ also misses to include the remaining $X_{i}$%
's corresponding to $i\in [m_{n}\ell _{n}+1, m_{n}\ell _{n}+r_{n}]$, where $m_{n}$ is the
integer part of $n/\ell _{n}$ and $r_{n}=n-m_{n}\ell _{n}.$ And although it
is possible to get rid of the term $\ell Y_{m+1,\ell }^{\ast}$ in line -1 of
their page 108 as we explained in the lines following the remark $(R2)$ in
the proof of Theorem \ref{theo1} below, we still think it would be handled
in the proof of Formula (\ref{C2}), as we did at the stage of Formula (\ref{C1A})
of the same proof below.\\

\noindent Based on this remark, the hypotheses (\ref{pauloB1}) and (\ref{pauloB2}) are
true and our ($Hab)$ holds. Then Formula (\ref{pauloB2S}) holds and Formula (\ref{C2}) is true. The Feller-Levy theorem handles the remaining part. Further%
\begin{equation*}
\frac{\ell (n)^{2}}{s_{n}^{2}}Var\left( Y_{m(n)+1,\ell (n)}^{\ast }\right)
\leq \frac{1}{s_{n}^{2}}\sum_{i=m(n)\ell (n)+1}^{r(n)}var(X_{i})+\frac{2\ell
(n)}{s_{n^{2}}}\left\{ \frac{1}{\ell (n)}\sum_{i=1}^{\ell (n)}u(i)\right\} .
\end{equation*}

\noindent Then, if Assumption (\ref{pauloB1}) and (\ref{pauloB2S}) hold, then ($Hab)$ is
implied by a general condition of the form

\begin{equation}
\frac{1}{s_{n}^{2}}\sum_{i=t_{n}}^{u_{n}}var(X_{i})\rightarrow 0\text{ as }%
n\rightarrow \infty ,  \label{HNab}
\end{equation}%
for $0\leq t_{n}\leq u_{n}\leq n,u_{n}-t_{n}\leq \ell
(n),(u_{n}-v_{n})/n\rightarrow 0$ as $n\rightarrow \infty$.\\

\bigskip \noindent \textbf{(2) With Cox-Grimmet Theorem \ref{cox}}.\\

\noindent It is immeadiate that the first part of Assumption (\ref{coxA1}) in that
theorem, that is
\begin{equation*}
Var(X_{j})\geq c_{1}>0,
\end{equation*}

\noindent implies, by association, that 
\begin{equation*}
s_{n}^{2}\geq \sum_{i=1}^{n}var(X_{i})\geq nc_{1}
\end{equation*}

\noindent and (\ref{pauloB2S}) holds since

\begin{equation*}
\lim \sup_{n\rightarrow +\infty }\frac{\ell (n)}{s_{n^{2}}}=\lim
\sup_{n\rightarrow +\infty }\frac{n}{s_{n^{2}}}\times \frac{\ell (n)}{n}\leq
c_{1}\lim \sup_{n\rightarrow +\infty }\frac{\ell (n)}{n}\leq c_{1}.
\end{equation*}

\noindent Next, the second part, that is 
\begin{equation*}
E\left\vert X_{j}\right\vert ^{3}\leq c_{2}<+\infty ,j\geq 1,
\end{equation*}

\noindent implies, by the formula $\left\vert x\right\vert ^{p}\leq 1+\left\vert
x\right\vert ^{q}$ for $1\leq p\leq q$ (see \cite{loeve}, page \ 157), that
for $c_{3}=1+c_{2}$,
\begin{equation*}
E\left\vert X_{j}\right\vert ^{3}\leq c_{3},j\geq 1,
\end{equation*}

\noindent and then Formula (\ref{HNab}) above holds since%
\begin{equation*}
\frac{1}{s_{n}^{2}}\sum_{i=t_{n}}^{u_{n}}var(X_{i})\leq c_{3}\frac{%
(u_{n}-t_{n})}{s_{n}^{2}}\leq \frac{c_{3}}{c_{1}}\frac{(u_{n}-t_{n})}{n}%
\rightarrow 0.
\end{equation*}

\noindent Next, by re-making the considerations given in Subsubsection \ref{subsubsec422}, Formula (\ref{fellerlevy01c}) of Theorem \ref{theo2} holds if 
\begin{equation*}
\frac{1}{s_{n}^{2}}\sum_{j=1}^{m(n)}\int_{\left\{ \left\vert
X_{j}\right\vert \geq \varepsilon s_{n}\right\} }X_{j}^{2}d\mathbb{P}%
\rightarrow 0\text{ as }n\rightarrow +\infty ,
\end{equation*}

\noindent for any $\varepsilon >0$. But we have under Condition \ref{coxA1} \ of
Theorem \ref{cox},

\begin{eqnarray*}
\frac{1}{s_{n}^{2}}\sum_{j=1}^{m(n)}\int_{\left\{ \left\vert
X_{j}\right\vert \geq \varepsilon s_{n}\right\} }X_{j}^{2}d\mathbb{P} &%
\mathbb{=}&\frac{1}{s_{n}^{2}}\sum_{j=1}^{m(n)}\int_{\left\{ \left\vert
X_{j}\right\vert \geq \varepsilon s_{n}\right\} }\frac{\left\vert
X_{j}\right\vert ^{3}}{\left\vert X\right\vert }d\mathbb{P} \\
&\leq &\frac{1}{\varepsilon s_{n}^{3}}\sum_{j=1}^{m(n)}\int_{\left\{
\left\vert X_{j}\right\vert \geq \varepsilon s_{n}\right\} }\left\vert
X_{j}\right\vert ^{3} \\
&\leq &\frac{m(n)c_{2}}{\varepsilon s_{n}^{3}}\leq \frac{c_{2}}{c_{1}^{3/2}}%
\frac{m(n)}{n^{3/2}}\rightarrow 0.
\end{eqnarray*}

\noindent Hence, Condition (\ref{fellerlevy01c}) is true. Finally Condition (\ref{coxA2})
ensures (\ref{C2}) and the Cox Theorem \ref{cox} is obtained.

\subsubsection{Conclusion} \label{subsubsec424} We conclude in two points.\\

\noindent \textbf{(A)} By combining our results with especially those of Oliveira \textit{et al.} \cite{paulo}, we have proved that the Newman approach already gave the best results in a Lyapounov-Feller-Levy type of \textit{CLT}. It is still possible to find different more or less sharp expressions of Conditions (\ref{C1}) and (\ref{C2}), stated in Subsection \ref{subsec42}. But no very different results cannot
be expected there. Yet, the \textit{CLT} problem is largement open since the current
results use the Newman approach. Is it possible to get rid of this approach
and and to use another one more general to establish more general \textbf{CLT}'s?
This seems to be the direction to be taken.\\ 

\noindent \textbf{(B)} In \cite{LAH2016}, an associated sequence is studied as a particluar case. Using a
direct method based on the characteristic function methid, it has been shown to
satistify the \textit{CLT} property. Yet, this sequence did not satisfy the
Cox-Grimmet condition $\inf_{n\geq 1}\mathbb{E}X_{i}^{2}\geq c_{1}>0$. This
kind of work may constitute a lead to more general \textit{CLT}.\\

\subsection{Proof of Theorem \ref{theo1}} \label{subsec5} As almost all the proofs of \textit{CLT}'s for associated
or weakly associated rv's, our proof is based on the three steps of the
original method of Newman and Wright (see \cite{newmanwright}). For compact
notation sake, we simply set $\ell (n)=\ell $ and $m(n)=m$. Let us define $%
\Psi _{\frac{S_{n}}{s_{n}}}(t)=\mathbb{E}\left( e^{itS_{n}/s_{n}}\right) ,$ $%
t\in \mathbb{R}$.\\

\noindent First, we have for $t\in \mathbb{R}$, 
\begin{equation*}
\left\vert \Psi _{_{\frac{S_{n}}{s_{n}}}}(t)-\Psi _{_{\frac{S_{m\ell }}{s_{n}%
}}}(t)\right\vert =\left\vert \mathbb{E}(e^{itS_{n}/s_{n}})-\mathbb{E}%
(e^{itS_{m\ell }/s_{n}})\right\vert
\end{equation*}
\begin{equation*}
=\left\vert \mathbb{E}\left[ e^{itS_{m\ell }/\sqrt{m\ell }}\left( e^{it\left[
(S_{n}/s_{n})-(S_{m\ell }/s_{n})\right] }-1\right) \right] \right\vert
\end{equation*}

\begin{equation}
\leq \mathbb{E}\left\vert e^{it\left( \frac{S_{n}}{s_{n}}-\frac{S_{m\ell }}{%
s_{n}}\right) }-1\right\vert .  \label{b}
\end{equation}

\noindent But for any $x\in \mathbb{R}$, 
\begin{equation*}
\left\vert e^{ix}-1\right\vert =|(\cos x-1)+i\sin x|=\left\vert 2\sin \frac{x%
}{2}\right\vert \leq |x|.
\end{equation*}

\bigskip \noindent Thus the second member of $(\ref{b})$ is, by the
Cauchy-Schwarz's inequality, bounded by 
\begin{equation*}
|t|\mathbb{E}\left\vert \frac{S_{n}}{s_{n}}-\frac{S_{m\ell }}{s_{n}}%
\right\vert \leq |t|\mathbb{V}ar\left( \frac{S_{n}}{s_{n}}-\frac{S_{m\ell }}{%
s_{n}}\right) ^{\frac{1}{2}}
\end{equation*}

\noindent and 
\begin{equation*}
\delta _{m,\ell }=\mathbb{V}ar\left( \frac{S_{n}}{s_{n}}-\frac{S_{m\ell }}{%
s_{n}}\right) =\frac{1}{s_{n}^{2}}\mathbb{V}ar\left( S_{n}-S_{m\ell }\right),
\end{equation*}

\noindent which tends to zero as $n\rightarrow +\infty $  by $(Hb)$ since%

\begin{eqnarray}
\delta _{m,\ell } &=&\frac{1}{s_{n}^{2}}\mathbb{V}ar\left( \
\sum\limits_{i=1}^{r}X_{m\ell +i}\right)   \label{C1A} \\
&\leq &\frac{\ell }{s_{n}^{2}}\mathbb{V}ar\left( \frac{1}{\sqrt{\ell }}\
\sum\limits_{i=1}^{\ell }X_{m\ell +i}\right)  \\
&\leq &C_{1}(n)\rightarrow 0.
\end{eqnarray}

\noindent This proves that 
\begin{equation}
|\Psi _{\frac{S_{n}}{s_{n}}}(t)-\Psi _{\frac{S_{m\ell }}{s_{n}}%
}(t)|\rightarrow 0\text{ as }n\rightarrow +\infty .  \label{commonStep01}
\end{equation}

\noindent \textbf{(R1)} Remark also for the purpose of Theorem \ref{theo2} that the
same conclusion holds when $(Hab)$ is true and we do not need (\textit{Hb}) in addition.\\

\noindent Next, remind that  $Y_{j,\ell }=(S_{j\ell }-S_{\ell (j-1)})/\sqrt{\ell }$%
, for $1\leq j\leq m$. \ Observe that%
\begin{equation*}
\frac{S_{m\ell }}{s_{n}}=\frac{\sqrt{\ell }}{s_{n}}\sum_{j=1}^{m}Y_{j,\ell }.
\end{equation*}

\noindent According to the Newman's inequality (see Lemma \ref{lemg3}), we have
\begin{equation*}
\left\vert \Psi _{\frac{S_{m\ell }}{s_{n}}}(t)-\prod\limits_{j=1}^{m}\Psi
_{Y_{j,n}}\left( \frac{\sqrt{\ell }}{s_{n}}t\right) \right\vert \leq \frac{%
\ell t^{2}}{2s_{n}^{2}}\sum_{1\leq j\neq k\leq m}Cov(Y_{j,\ell },Y_{k,\ell
}).
\end{equation*}

\noindent But, 
\begin{eqnarray*}
\frac{\ell t^{2}}{2s_{n}^{2}}\sum_{1\leq j\neq k\leq m}Cov(Y_{j,\ell
},Y_{k,\ell }) &=&\frac{\ell t^{2}}{2s_{n}^{2}}\mathbb{V}ar\left(
\sum_{j=1}^{m}Y_{j,\ell }\right) -\frac{\ell t^{2}}{2s_{n}^{2}}\sum_{j=1}^{m}%
\mathbb{V}ar(Y_{j,\ell }) \\
&=&\frac{t^{2}}{2}\left[ \mathbb{V}ar\left( \frac{\sqrt{\ell }}{s_{n}}%
\sum_{j=1}^{m}Y_{j,\ell }\right) -\frac{\ell }{s_{n}^{2}}\sum_{j=1}^{m}%
\mathbb{V}ar\left( Y_{j,\ell }\right) \right]  \\
&=&\frac{t^{2}}{2}\left[ \mathbb{V}ar\left( \frac{1}{s_{n}}S_{m\ell }\right)
-\frac{\ell }{s_{n}^{2}}\sum_{j=1}^{m}\mathbb{V}ar\left( \frac{S_{j\ell
}-S_{\ell (j-1)}}{\sqrt{\ell }}\right) \right]  \\
&\leq &\frac{t^{2}}{2}\left[ 1-\frac{\ell }{s_{n}^{2}}\sum_{j=1}^{m}\mathbb{V%
}ar\left( \frac{S_{j\ell }-S_{\ell (j-1)}}{\sqrt{\ell }}\right) \right]\\
&-&\frac{t^{2}}{2s_{n}^{2}}\mathbb{V}ar\left( \ \sum\limits_{j=m\ell
+1}^{n}X_{j}\right) ,
\end{eqnarray*}

\noindent which tends to zero as $n\rightarrow +\infty $ by \textit{(Ha)} and
\textit{(Hb)}, that is 
\begin{equation}
\left\vert \Psi _{\frac{S_{m\ell }}{s_{n}}}(t)-\prod\limits_{j=1}^{m}\Psi
_{Y_{j,n}}\left( \frac{\sqrt{\ell }}{s_{n}}t\right) \right\vert \rightarrow 0%
\text{ as }n\rightarrow +\infty .  \label{conclusion01}
\end{equation}

\noindent The proof will be completed by establishing that%
\begin{equation}
\prod\limits_{j=1}^{m}\Psi _{Y_{j,n}}\left( \frac{\sqrt{\ell }}{s_{n}}%
t\right) \rightarrow \exp (-t^{2}/2)\text{ as }n\rightarrow +\infty .
\label{lastStep}
\end{equation}

\noindent \textbf{(R2)} Here, we make a second remark which is relevant to the Proof
of Theorem \ref{theo2} and next to generalisations of the results. The above
computations led to

\begin{eqnarray*}
0 &\leq &\frac{\ell t^{2}}{2s_{n}^{2}}\sum_{1\leq j\neq k\leq
m}Cov(Y_{j,\ell },Y_{k,\ell })=\frac{t^{2}}{2}\left[ 1-\frac{\ell }{s_{n}^{2}%
}\sum_{j=1}^{m}\mathbb{V}ar\left( \frac{S_{j\ell }-S_{\ell (j-1)}}{\sqrt{\ell }}\right) \right]\\
&-&\frac{t^{2}}{2s_{n}^{2}}\mathbb{V}ar\left( \
\sum\limits_{j=m\ell +1}^{n}X_{j}\right)  \\
&\leq &\frac{t^{2}}{2}\left[ 1-\frac{\ell }{s_{n}^{2}}\sum_{j=1}^{m}\mathbb{V%
}ar\left( \frac{S_{j\ell }-S_{\ell (j-1)}}{\sqrt{\ell }}\right) \right] .
\end{eqnarray*}

\noindent Then only \textit{(Ha)} is needed to ensure \ref{lastStep}.\\

\noindent Now, we resume to the normal course of our demonstration. From this step,
the conclusion on the weak law of $S_{n}/s_{n}$, comes uniquely from Formula
(\ref{lastStep}) which  expresses the weak convergence
of sums of the form
\begin{equation}
T_{m(n)}^{\ast }=\frac{1}{s_{n}}\sum_{j=1}^{m(n)}V_{j}, \label{sumIndep}
\end{equation}

\noindent where the the $V_{j}$'s are independent random variables such that
for each $j\in \{1,m\},V_{j}^{\ast }$ has the same law as $S_{j\ell
}-S_{(j-1)\ell }. $ Remind that, for each $j\in \{1,m\},\tau
_{j}^{2}=Var\left( S_{j\ell }-S_{(j-1)\ell }\right) =\mathbb{E}\left(
S_{j\ell }-S_{(j-1)\ell }\right) ^{2}$ and 
\begin{equation*}
\nu _{m(n)}^{2}=\tau _{1}+...+\tau _{m(n)},n\geq 1\text{.}
\end{equation*}

\noindent By Assumption $(Ha),$ we have $\nu _{m(n)}/s_{n}\rightarrow 1$ as $%
n\rightarrow +\infty $ and by Slustsky theorem (see for example Proposition
15 in \cite{wcia-srv-ang}, page 60), the weak convergence, if it holds,
would be the same as that of

\begin{equation*}
T_{m(n)}=\frac{1}{v_{m(n)}}\sum_{j=1}^{m(n)}V_{j}.
\end{equation*}

\noindent Condition (\textit{Hb}) is the Lyapounov's one for this problem (see Lo\`{e}ve \cite%
{loeve}, page 287, Point B), where $v_{m(n)}$ is replaced by $s_{n}$.
\noindent This completes the proof.\newline

\subsection{Proof of Theorem \ref{theo2}}

Based of the remarks marqued (\textbf{R1}) and \textbf{(R2)} in the body of the proof of Theorem %
\ref{theo1}, we conclude that if \textit{(L)}, \textit{(Ha) and (Hab) hold, }%
the conclusion on the weak law of $S_{n}/s_{n}$, comes uniquely from
Formula (\ref{lastStep}). At this step, the condition on the $(2+\delta )^{th}$
moments, that $\mathbb{E}\left\vert X_{j}\right\vert ^{2+\delta }<+\infty $,$%
j\geq 1,$ is not required. And, Formula \ref{lastStep} expresses the weak
convergence of the sums defined in (\ref{sumIndep}).\\

\noindent From there, the problem becomes the classical
Lyapounov-Levy-Feller Theorem. And we have the following conclusion :

\noindent \textbf{(a)} $\max_{1\leq k\leq m(n)}\{\tau _{j}/\nu
_{m(n)}\}\rightarrow 0$ as $n\rightarrow +\infty $ and 
\begin{equation*}
\frac{1}{v_{m(n)}}\sum_{j=1}^{m(n)}V_{j}\rightsquigarrow \mathcal{N}(0,1)%
\text{ as }n\rightarrow +\infty ,
\end{equation*}

\noindent if and only if \newline

\noindent \textbf{(b)} for any $\varepsilon >0,$ 
\begin{equation*}
g(\varepsilon )=\frac{1}{v_{m(n)}^{2}}\int_{(\left\vert x\right\vert \geq
\varepsilon v_{m(n)})}x^{2}dF_{V_{j}}\rightarrow 0\text{ as }n\rightarrow
+\infty .
\end{equation*}

\bigskip \noindent These two conditions are exactly those given in the
statement of the theorem, where the replacement of $(\left\vert x\right\vert
\geq \varepsilon v_{m(n)})$ by $(\left\vert x\right\vert \geq \varepsilon
s_{n})$ in the expression of $g$ is possible because of $\nu
_{m(n)}/s_{n}\rightarrow 1$ as $n\rightarrow +\infty $.\newline

\noindent This finishes the proof on this theorem.

\subsection{Proof of Theorem \protect\ref{theo3}}

The proof of Theorem \ref{theo3} is based on that of Theorem \ref{theo1}
from Formula (\ref{conclusion01}).\\

\textbf{Acknowledgment} The second author acknowledges support from the
World Bank Excellence Center (CEA-MITIC) that is continuously funding his
research activities starting 2014. The first author thanks the \textbf{%
Programme de formation des formateurs} of USSTB who financed his stays in
the LERSTAD of UGB while preparing his Ph.D dissertation. Both authors
acknowledge support from the \textbf{R\'{e}seau EDP - Mod\'{e}lisation et
Contr\^{o}le}, of Western African Universities, that financed travel and
accomodation of the second author while visiting USTTB in preparation of a
series of works with his PhD students there.


\begin{thebibliography}{99}
\bibitem{bulinski2007} Bulinski A. and Shashkin A.(2007). Limit theorems for
associated random fields and related systems. World Scientific Publishing,
Singapore.

\bibitem{burton} Burton, R.M., Dabrowski, A.R. and Dehling, H. (1986). An
invariance principle for weakly associated random variables, \textit{Stoch.
Proc. Appl.}, 23, 301-306.

\bibitem{cox} Cox, J.T. and Grimmett, G. (1984) Central limit theorems for
associated random variables and the percolation model, \textit{Ann. Probab.}%
, 12, 514-528.

\bibitem{dabro} Dabrowski, A.R. and Dehling, H. (1988). A Berry-Esseen
theorem and a functional law of the iterated logarithm for weakly associated
random variables, \textit{Stochastic Process. Appl.}, 30, 247-289.

\bibitem{daley} Daley, D. J.(1968). Stochastically monotone Markov Chains,
Z. \textit{Wahrsch. theor. verw Gebiete}, 10, 305-317.

\bibitem{esary} Esary, J., Proschan, F. and Walkup, D.(1967). Association of
random variables with application. \textit{Ann. Math Statist.}, 38,

\bibitem{fortuin} Fortuin, C., Kastelyn, P. and Ginibre, J.(1971).
Correlation inequalities on some partially ordered sets. \textit{Comm. Math.
Phys.}, 22, 89-103.

\bibitem{gut} Gut, A. (2005). \textit{Probability : A Graduate Course}.
Springer Science+Business Media, Inc. ISBN 0-387-22833-0.

\bibitem{jiazhu} Jiazhu, P.(2002). Tail dependence of random variables from
ARCH and heavy-tailed bilinear models. \textit{Sciences in China}, 45 (6),
Ser. A, 749-760$.$

\bibitem{cam} Le Cam, L.(1986). The Central Limit Theorem Around 1935. \textit{Statistical Science.}, Vol. 1 (1), 78-91.

\bibitem{LAH2016} Lo G.S., Fall A.M. and Harouna S.(2016). A Central limit Theorem
of dependent sums of standard exponential functionals motivated by extreme
value theory. \textit{To come}.

\bibitem{lehmann} Lehmann, E. L.(1966). Some Concepts of dependence. \textit{%
\ Ann. Math. Statist.}, 37, 1137-1153.

\bibitem{wcia-srv-ang} Lo, G.S.(2016). Weak Convergence (IA). Sequences of
random vectors. SPAS Books Series.(2016). Doi : 10.16929/sbs/2016.0001.
Arxiv : 1610.05415

\bibitem{loeve} Lo\`{e}ve, M.(1977). \textit{Probability Theory I}.
Springer-Verlag. New-York.

\bibitem{louhichi} Louhichi, S.(2000). Weak convergence for empirical
processes of associated sequences. Ann. Inst. Henri Poincar\'e,
Probabilit\'es et Statistiques 36 (\textbf{5}), pp. 547\^{a}\euro ``567.

\bibitem{newman80} Newman C.M. (1980) Normal fluctuations and the FKG
inequalities. \textit{Comm. Math. Phys.} 74, 119-128.

\bibitem{newmanwright} Newman, C.M and Wright, A.L.(1981). An invariance
principle for certain dependent sequences. \textit{Ann. probab.}, \textbf{9}
(4), 671-675.

\bibitem{newmanwright2} Newman, C.M and Wright, A.L.(1982). Associated
random variables and martingale inequalities. Z. Wahrscheinlichkeitstheorie
verw. Gebiete 59, 361-371.

\bibitem{paulo} Oliveira, P.E.(2012). \textit{Asymptotics for Associated
Random Variables.} DOI 10.1007/978-3-642-25532-8,\textit{\ }\copyright\ %
Springer-Verlag Berlin Heidelberg.

\bibitem{pitt} Pitt, L.(1982). Positively Correlated normal variables are
associated. \textit{Ann. Probab.}, 10, 496-499.

\bibitem{rao} Prakasa Rao, B. L. S.(2012). \textit{Associated sequences,
Demimartingales and Nonparametric Inference. Probability and its
applications }. Springer Basel Doredrecht, Heidelberg, London, New York.

\bibitem{gslo} Sangar\'{e}, H. and Lo, G. S. A Review on asymptotic
normality of sums of associated random variables. \textit{Afrika Statistika}%
, 11 (\textbf{1}), pp.855-867. Doi : 10.16929/as/2016.855.79. Arxiv
1405.4316.

\bibitem{yu93} Yu, H.(1993). A Gkivenko-Cantelli lemma and weak convergence
for empirical processes of associated sequences. \textit{Probab. Theory
Related Fields.} 95, 357-370.
\end{thebibliography}
\end{document}